\documentstyle [12pt,amssymb,amscd,psfig]{amsart}

\newtheorem{theorem}{Theorem}[section]
\newtheorem{lemma}[theorem]{Lemma}
\newtheorem{corollary}[theorem]{Corollary}

\theoremstyle{definition}
\newtheorem{definition}[theorem]{Definition}

\theoremstyle{remark}
\newtheorem{remark}[theorem]{Remark}

\numberwithin{equation}{section}

\begin{document}

\title[$2$-Complexes in ${\mathbb{R}}^{4}$]
{Embedding Obstructions and {\boldmath $4$}-dimensional Thickenings
of {\boldmath $2$}-complexes}

\author{Vyacheslav S. Krushkal}
\address{Department of Mathematics, Yale University, New Haven CT 06520}
\email{krushkal\char 64 math.yale.edu} 

\begin{abstract}
The vanishing of Van Kampen's obstruction is known to be necessary
and sufficient for embeddability of a simplicial $n$-complex into
${\mathbb{R}}^{2n}$ for $n\neq 2$, and it was recently shown to be incomplete
for $n=2$. We use algebraic-topological invariants of four-manifolds
with boundary to introduce a sequence of higher embedding obstructions
for a class of $2$-complexes in ${\mathbb{R}}^4$.
\end{abstract} 


\maketitle

\section{Introduction} \label{introduction} 

By general position any $n$-dimensional simplicial complex $K$
PL embeds into ${\mathbb{R}}^{2n+1}$, while
the image of a generic map of $K$ into ${\mathbb{R}}^{2n}$ has
a finite number of double points.
By counting double points of an immersion one gets
the cohomological obstruction to embeddability of an $n$-complex
into ${\mathbb{R}}^{2n}$, introduced by Van Kampen \cite{Van Kampen}.
He also constructed for each $n$ examples which do not admit an embedding.
An application of Whitney trick shows that this obstruction
is complete for $n>2$, see \cite{Van Kampen}, 
\cite{Shapiro}, \cite{Wu}, \cite{FKT}.
It follows from Kuratowski's planarity criterion for graphs
\cite{Kuratowski} that this result also holds for $n=1$.
The remaining case, $n=2$, was open until recently when the
obstruction was shown in \cite{FKT} to be incomplete. 

This paper is centered around the question of
embeddability of $2$-complexes in ${\mathbb{R}}^4$, and
is motivated by the result of \cite{FKT}.
We define for $2$-complexes $K$ with $H_1(K;{\mathbb{Q}})=0$ a sequence 
of higher embedding obstructions $\{o_m(K)\}$, using Massey 
products on the boundary of a $4$-dimensional thickening $M^4$ of $K$. 
Roughly, Van Kampen's obstruction 
corresponds in this setting to the intersection pairing on $M$, 
modulo the choice of a thickening $M$. 
Since different thickenings may give different Massey products,
$\{o_m(K)\}$ in general are {\em subsets} of the corresponding
cohomology groups; $o_{m+1}(K)$ is defined if $o_m(K)$
contains zero. If $K$ embeds into ${\mathbb{R}}^4$ then 
$0\in o_m(K)$ for each $m$.
We prove that these higher obstruction detect non-embeddability of 
the family of examples introduced in \cite{FKT}, by showing that $o_m(K)$
does not contain zero for some $m$.
Our proof uses the result of Conway - Gordon and Sachs that any embedding
of a complete graph on $6$ vertices into $S^3$ contains two disjoint
linking cycles (\cite{CG}, \cite{Sachs}). 

In the simplest relative case, for the disjoint union of 
$2$-disks with a prescribed embedding of their boundaries
into $S^3$, by a result of Turaev \cite{Turaev} Massey products
correspond to Milnor's $\bar\mu$-invariants of the link in $S^3$, 
so our obstructions may be thought of as an absolute analogue
of $\bar\mu$-invariants.
As in the case of $\bar\mu$-invariants (for example, Whitehead 
double of the Hopf link is not a slice link, while all $\bar\mu$-invariants 
vanish), one does not expect that the entire sequence of obstructions
defined here is complete, although no examples are known at this time.
The question about $2$-complexes has an additional subtlety, being in  
piecewise-linear category, where embeddings are not necessarily locally 
flat.

The definition of Van Kampen's obstruction is recalled in section 
\ref{VKobstruction}. In section \ref{thickenings} we prove its 
reformulation in the context of thickenings, and we introduce the 
sequence of higher obstructions $\{o_m(K)\}$. We review the examples
of $2$-complexes in \cite{FKT} in section \ref{examples}, and we
compute the obstructions for these examples. Section \ref{config}
gives a reformulation of Van Kampen's obstruction in terms of
configuration spaces, which suggests another approach to defining
higher embedding obstructions.

The present study of the embedding problem for $2$-complexes in 
${\mathbb{R}}^{4}$ is motivated, in part, by the $4$-dimensional topological 
surgery conjecture, via its {\em (A,B)-slice}
reformulation \cite{Freedman2}. More precisely, the surgery conjecture
is equivalent to the relative embedding question for a certain family of 
$4$-dimensional handlebodies -- ``thickenings'' of $2$-complexes
in the sense of section \ref{thickenings}. However, many interesting 
examples of these handlebodies have non-trivial first homology, 
and for this application the obstructions $\{o_i(K)\}$ 
need to be extended to the general case.

\section{Van Kampen's obstruction} \label{VKobstruction} 

In this section we briefly review the definition of Van Kampen's
obstruction, more details are given in \cite{FKT}.
In 1933 Van Kampen \cite{Van Kampen} introduced an obstruction
$o(K)\in H^{2n}_{{{\mathbb{Z}}}/{2}}(K^*;{\mathbb{Z}})$ to piecewise-linear 
embeddability
of an $n$-dimensional simplicial complex $K$ into ${{\mathbb{R}}}^{2n}$.
The cohomology in question is ${{\mathbb{Z}}}/2$-equivariant cohomology
where ${{\mathbb{Z}}}/2$ acts on the deleted product 
$K^*\!=K\times K\setminus\Delta$ of a complex $K$
by exchanging the factors of $K^*$ and acts on the coefficients
by multiplication with $(-1)^n$.  
The diagonal $\Delta$ consists of all products $\sigma\times\tau$
such that simplices $\sigma$ and $\tau$ have at least one vertex in common.
Note that for $n$ even (in particular, in the case of main interest
in this paper, $n=2$) the action of ${\mathbb{Z}}/2$ on the coefficients 
is trivial, and $o(K)$ is an element of the ordinary 
cohomology group $H^{2n}(K^{*}/({\mathbb{Z}}/2); {\mathbb{Z}})$.

Let $f$ be any PL immersion of $K$ into ${{\mathbb{R}}}^{2n}$. The obstruction 
is defined on the cochain level by counting 
algebraic intersection numbers of the images of disjoint top-dimensional
simplices of $K$: $o_f(\sigma^n\times\tau^n)=f(\sigma)\cdot f(\tau)$.
Here $\sigma\times\tau$ is viewed as an oriented generator of the
$2n$-th chain group of $K\times K\setminus\Delta$. 
The cohomology
class $o(K)$ of $o_f$ is independent of the chosen immersion $f$.
Clearly $o(K)$ is trivial if $K$ embeds into ${{\mathbb{R}}}^{2n}$.
Shapiro \cite{Shapiro} and Wu \cite{Wu} made this definition precise
and proved, using the Whitney trick, the converse in dimension $n$
greater than $2$.

\begin{theorem}[\cite{Van Kampen}, \cite{Shapiro}, \cite{Wu}, 
\cite{Kuratowski}] \label{iff}
For $n\neq 2$ an $n$-dimensional simplicial complex $K$ admits an  
embedding into ${{\mathbb{R}}}^{2n}$ if and only if $o(K)$ vanishes.
\end{theorem}

See \cite{FKT} for a modern exposition of the proof for $n>2$. 
For $n=1$ this theorem follows from Kuratowski's planarity criterion.
The obstruction in the remaining case, for $n=2$, was shown to be
incomplete in \cite{FKT}. We recall the construction
of examples in \cite{FKT} in section \ref{examples}.

\section{Obstructions via $4$-dimensional thickenings
of $2$-complexes} \label{thickenings}  

In this section we give a rational reformulation of Van Kampen's 
obstruction $o(K)$ in terms of thickenings of $K$, and we introduce 
a sequence of higher embedding obstructions $\{o_m(K)\}$ for 
$2$-complexes whose rational first homology vanishes.
Throughout this section all coefficients are ${\mathbb{Q}}$,
unless stated otherwise, and $K$ denotes a simplicial $2$-complex.

\begin{definition} \label{thickening}
A {\it thickening} of $K$ is 
a smooth $4$-manifold $M$ with boundary, obtained by replacing
each $i$-simplex of $K$ with a $4$-dimensional $i$-handle, $i=0,1,2$.
The attaching map of each $2$-handle is
required to be isotopic, within the union of $0$- and $1$-handles,
to the attaching map of the corresponding $2$-dimensional simplex.
\end{definition}

In general, $K$ may have different thickenings depending on the choice of 
attaching maps of the $2$-handles. For example, $S^2\times S^2\setminus 
4$-cell and the boundary-connected sum $S^2\times D^2\natural 
S^2\times D^2$ are both thickenings of $S^2\vee S^2$. 

The intersection pairing on $M$ defines an element 
$\iota\!\in Hom(H_2(M)\otimes H_2(M),\mathbb{Q})$. 
Let $\bar\iota\in H^4(K\times K\setminus{\Delta};{\mathbb{Q}})$ denote 
its image under the homomorphism 

\[ Hom(H_2(M)\otimes H_2(M),{\mathbb{Q}})\cong 
Hom(H_2(K)\otimes H_2(K),{\mathbb{Q}})\cong \]

\[ \cong H^4(K\times K) 
\longrightarrow H^4(K\times K\setminus{\Delta} ) \]

\noindent
where the last map is induced by inclusion.

\begin{theorem} \label{reformulation2}  \sl
The image of the (rational) Van Kampen's obstruction $o(K)$
under the homomorphism induced by the quotient map

\[ H^{4}_{{\mathbb{Z}}/{2}}(K\times K\setminus{\Delta};{\mathbb{Q}})
\longrightarrow H^4(K\times K\setminus{\Delta};{\mathbb{Q}}) \]

\noindent
coincides with $-\bar\iota$.
\end{theorem} 

\noindent
{\em Proof.} Suppose a thickening $M$ is induced by an immersion
$f: K\longrightarrow {\mathbb{R}}^4$, so that $f$ extends to an immersion
$M\longrightarrow {\mathbb{R}}^4$. 
By subdividing the complex $K$, if necessary, one may assume
that $f(\sigma)\cap f(\tau)= \emptyset$ for all (open) $2$-simplices
$\sigma\neq\tau$ with $\sigma\times\tau\in\Delta$,
and $f|_{\sigma}$ is an embedding for each $\sigma$.
Let $\bar o_f: C_4(K\times K)\longrightarrow \mathbb{Q}$
denote the extension by zero on the diagonal of Van Kampen's
cochain $o_f: C_4(K\times K\setminus\Delta)\longrightarrow \mathbb{Q}$. 
It suffices to prove that $[\bar o_f ]$
and $-\iota$ define identical elements in 
$Hom(H_2(K)\otimes H_2(K),\mathbb{Q})$. Let $a$, $b$ be two classes 
in $H_2(K)$ and let $\alpha=\Sigma\alpha_i\sigma_i$,
$\beta=\Sigma\beta_i\sigma_i$ be their cycle representatives,
where $\{\sigma_i\}$ is the set of $2$-simplices of $K$.
In order to compute $a\cdot b$ in $M$, perturb $\alpha$
and $\beta$ to $\tilde{\alpha}$ and $\tilde{\beta}$ which
intersect each other transversely (in a finite number of
double points). The intersection number of two cycles 
$f(\tilde{\alpha})$ and $f(\tilde{\beta})$ in ${\mathbb{R}}^4$ is
trivial. On the other hand, $f(\tilde{\alpha})\cdot f(\tilde{\beta})$
may be computed as the sum of two terms: one is the intersection number of 
$\tilde{\alpha}$ and $\tilde{\beta}$ in $M$, the other is obtained by
considering the intersections of $f(\tilde{\alpha})$ and $f(\tilde{\beta})$
in ${\mathbb{R}}^4$,
which are singular points of $f$. This last term is equal to
$\bar o_f(\alpha\times\beta)$, and this proves 

\[ [\bar o_f]=-\iota: H_2(K)\otimes H_2(K)\longrightarrow\mathbb{Q}. \]

\noindent
The restriction of $[\bar o_f ]$ to 
$H^{4}(K\times K\setminus\Delta; \mathbb{Q})$
coincides with $o(K)$, thus the result is proved for thickenings
induced by immersions.

In general not every thickening of $K$ may be immersed into
${\mathbb{R}}^4$. Let $M^4$ be an arbitrary thickening of $K$
and let $f: K\longrightarrow {\mathbb{R}}^4$ be any immersion.
Again one may assume that $f(\sigma)\cap f(\tau)=\emptyset$
if $\sigma\times\tau\in\Delta$, $\sigma\neq\tau$, and $f|_{\sigma}$
is an embedding for each simplex $\sigma$.
The immersion $f$ extends to an embedding of $0$- and $1$-handles of $M$.
There is an integer obstruction to extending it over each $2$-handle,
due to a possible difference of the framing of the $2$-handle
and of the normal bundle of the $2$-simplex in ${\mathbb{R}}^4$.
However, each $2$-handle may be mapped into ${\mathbb{R}}^4$
as a bundle over the corresponding $2$-simplex, pinched 
over several points.

The proof, given above in the case of an immersion, carries through,
if one extends $o_f$ to $\bar o_f$ by setting
$\bar o_f(\sigma\times\sigma)$ to be equal to the 
difference in framings, discussed above, and setting
$\bar o_f(\sigma\times\tau)=0$ for all
$\sigma\times\tau\in\Delta$, $\sigma\neq\tau$.
\qed

\begin{remark}
In general
the intersection pairing varies within the homotopy type of
a $4$-manifold $M$. In the example above
the intersection pairing on $S^2\times D^2\natural S^2\times D^2$
is trivial, while the pairing on $S^2\times S^2\setminus 4$-cell is
non-degenerate. However, theorem \ref{reformulation2} shows
that the pull-back of the intersection pairing on thickenings
to a cohomology class on 
$K\times K\setminus\Delta$ is an invariant of $K$,
which coincides with the image of the (negative) Van Kampen's obstruction.
\end{remark}

As a corollary to the proof of theorem \ref{reformulation2}, we have
the following result.

\begin{lemma} \label{zeropairing}  \sl
Let $K$ be a $2$-complex such that Van Kampen's obstruction $o(K)$ 
vanishes. Then there is 
a $4$-dimensional thickening $M$ of $K$ with the trivial 
intersection pairing $\iota=0\in Hom(H_2(M)\otimes H_2(M); \mathbb{Q})$.
\end{lemma} 

\noindent
{\em Proof.}
Any cochain representative of the obstruction $o(K)$ is given
by $o_f$ for some immersion $f$, see \cite{Van Kampen} or \cite{FKT}.
Since $o(K)$ vanishes, there exists an immersion $f: K\longrightarrow
{\mathbb{R}}^4$, giving rise to the trivial Van Kampen's cochain
$o_f=0$. Let $M$ denote the thickening induced by $f$.
It follows from the proof of theorem \ref{reformulation2}
that if one extends $o_f$ by zero on the diagonal to a
cochain $\bar o_f$ on $K\times K$, then $\iota=[\bar o_f]=0\in
Hom(H_2(M)\otimes H_2(M);\mathbb{Q})$.
\qed 

\vspace{.3cm}   

Before introducing the higher embedding obstructions, 
we recall the definition of Massey products. See 
\cite{Kraines} for proofs and additional properties.

\begin{definition}
Let $X$ be a space, and let $\alpha_1,\ldots,\alpha_m$ be elements in 
$H^1(X)$. Suppose there is a collection of 1-cochains
$S=\{c_{ij}\in C^1(X) | 1\leq i\leq j\leq m, (i,j)\neq (1,m)\}$
satisfying

\[ [c_{ii}]=\alpha_i\;\,\mathrm{for\;\,each}\;\,i=1,\ldots,m, \]

\[ \delta c_{ik}=\sum_{j=i}^{k-1} 
c_{ij}\cup c_{j+1,k}\;\,\mathrm{for}\;\,i<k. \]

\noindent
Then the cochain $\sum_{j=1}^{m-1} c_{1j}\cup c_{j+1,m}$
is a cocycle, and its cohomology class in $H^2(X)$ is called
the {\em Massey product} of $\alpha_1,\ldots,\alpha_m$ defined by 
the system $S$.
The set of Massey products corresponding to all such defining
systems is denoted by $<\!\alpha_1,\ldots,\alpha_m\!>\subset H^2(X)$.
\end{definition}

Massey product of two elements is just a cup product. 
Note that given some classes $\alpha_1,\ldots,\alpha_m$,
$<\!\alpha_1,\ldots,\alpha_m\!>$ is not necessarily defined.
However, if all Massey products of less than $m$ elements
vanish, then for any $\alpha_1,\ldots,\alpha_m\in H^1(X)$, 
$<\!\alpha_1,\ldots,\alpha_m\!>$ is a well-defined
element.
   
The following lemma justifies our definition of higher embedding obstructions.

\begin{lemma} \label{productsvanish}  \sl
Let $M$ be a $4$-manifold with boundary and with $H_1(M;{\mathbb{Q}})=0$,
and suppose $M$ admits an embedding into ${\mathbb{R}}^4$. 
Then all Massey products on $H^1(\partial M; \mathbb{Z})$ vanish.
\end{lemma}

\noindent
{\em Proof.} 
Let $N$ denote the complement ${\mathbb{R}}^4\setminus M$.
By Alexander duality, $H_2(N)$ and $H^2(N)$ are trivial.
The map $i^{*}: H^1(N)\longrightarrow H^1(\partial M)$
in the cohomology sequence of the pair $(N, \partial M)$
is an isomorphism, since by assumption and by Poincar\'{e} duality
$H^1(N, \partial M)\cong H_3(N)$ and $H^2(N, \partial M)\cong H_2(N)$
are trivial. Assume inductively that all Massey products of length
less than $m$ vanish for some $m\geq 2$; then for any 
${\alpha}_1,\ldots,{\alpha}_m\in H^1(\partial M)$ one has
 
\[ <\!{\alpha}_1,\ldots,{\alpha}_m\!>=
i^{*}\!\!<\!(i^{*})^{-1}{\alpha}_1,\ldots,
(i^{*})^{-1}{\alpha}_m\!>\,\in H^2(\partial M). \]

\noindent
However, this is the image of an element in $H^2(N)=0$,
and the result follows. \qed

\vspace{.3cm}

Let $K$ be a $2$-complex with $H_1(K; {\mathbb{Q}})=0$, and assume $o(K)$ 
vanishes. Let $M$ be a thickening of $K$ with 
trivial intersection pairing (its existence is given by
lemma \ref{zeropairing}.)
Note that the map $H_2(\partial M)\longrightarrow H_2(M)$
is an isomorphism, since by assumption on $K$, $H_3(M,\partial M)\cong
H^1(M)=0$, and the map $H_2(M)\longrightarrow H_2(M,\partial M)$
is trivial by assumption on the intersection pairing.

We now give the definition of higher embedding obstructions.
Let $a_1, a_2, a_3$ be classes in $H_2(K)$, and let 
$\alpha_1, \alpha_2, \alpha_3\in H^1(\partial M)$ denote
their images under the isomorphism 

\[ H_2(K)\cong H_2(M)\cong H_2(\partial M)\cong H^1(\partial M). \]

\noindent 
The triple cup product $(\alpha_1\cup\alpha_2\cup\alpha_3)[\partial M]$
defines a homomorphism $H_2(K)\otimes H_2(K)\otimes H_2(K)
\longrightarrow \mathbb{Q}$, and an element 
$o_3(K, M)\in H^6(K\times K\times K)$.
The cohomology class $o_3(K, M)$ depends in general on the choice of
a thickening $M$, thus we define the third obstruction $o_3(K)$ to be
the subset $\{o_3(K, M)\}\subset H^6(K^3; {\mathbb{Q}})$ where $M$ 
is to vary over all thickenings of $K$ with trivial intersection pairing.
Note that if $o_3(K)$ is defined and contains zero, then there is a 
thickening $M$ of $K$ such that all cup products on $H^2(\partial M)$ vanish.

\begin{definition}
Define $o_2(K)$ to be the Van Kampen's obstruction $o(K)$. If $o_2(K)$ 
vanishes, then $o_3(K)\subset H^6(K^3)$ is defined as above. Assume by 
induction that for some $m>3$ there is a thickening $M$ of $K$ such that
$o_{m-1}(K, M)$ is defined and is equal to zero (equivalently, the 
intersection pairing on $M$ is trivial, and all Massey products on 
$H^1(\partial M)$ of at most $(m-2)$ elements vanish.)
Let $a_1,\ldots, a_m$, be classes in $H_2(K)$, and
let $\alpha_1,\ldots, \alpha_m$ denote the corresponding
elements in $H^1(\partial M)$. 
The class $o_m(K,M)\in H^{2m}(K^m;\mathbb{Q})$
is defined by the homomorphism

\[ H_{2m}(K^m)\cong \otimes_1^m H_2(K)
\cong \otimes_1^m H^1(\partial M)
\longrightarrow \mathbb{Q} \]

\noindent
which sends $a_1\!\otimes\cdots\otimes\!a_m$ to 
$(<\!\alpha_1,\ldots,\alpha_{m-1}\!>\!\cup\,\alpha_m)[\partial M]$.
\end{definition}

Here since all Massey products on $H^1(\partial M)$ of less 
than $(m-\!1)$ elements vanish, $<\!\alpha_1,\ldots,\alpha_{m-1}\!>
\in H^2(\partial M)$ is a well-defined element.

\begin{definition}
The obstruction $o_m(K)$ is defined to be the subset

\[ \{o_m(K,M)\}\subset H^{2m}(K^m), \]

\noindent 
where $M$ is to vary over all thickenings such that $o_{m-1}(K,M)=0$.
Note that $o_m(K)$ is defined if $o_{m-1}(K)$ is defined and contains zero.
\end{definition}

Lemma \ref{productsvanish} implies the following corollary.

\begin{corollary} \sl
Let $K$ be a $2$-complex with $H_1(K;{\mathbb{Q}})=0$.
If $K$ admits an embedding into ${\mathbb{R}}^4$ then $o_m(K)$ is 
defined and contains zero for each $m$.
\end{corollary}

In section \ref{examples} we show that 
$o_m(K)$ does not contain zero for some $m$ for examples in \cite{FKT},
thus giving another proof that they do not embed into ${\mathbb{R}}^4$.


\noindent
{\bf The relative embedding problem.}
Let $K$ be a $2$-complex with $H_1(K;{\mathbb{Q}})=0$, and
let $L$ be a $1$-dimensional subcomplex of $K$ with a prescribed 
embedding $\phi: L\hookrightarrow S^3$. Consider the relative
embedding problem: does there exist an embedding $K\hookrightarrow B^4$
which extends $\phi$?
Denote $B^4\cup_{\phi}$(thickening of $K$) by $M$,
where thickening is taken in the sense of Definition \ref{thickening}. 
Let $K^m_{L}$ denote the subset in $K^m$ consisting of all
$m$-tuples $(x_1,\ldots,x_m)$ such that $x_i\in L$ for some $i$.
Assume that ${\pi}_0(L)\longrightarrow {\pi}_0(K)$ is injective
to have $H_1(M)=0$.
Analogously to the absolute case, Massey products on 
$\partial M$ define an element, depending on $M$, 
in the relative cohomology group $H^{2m}(K^m,K^m_{L})$. 
Let $o_m(K,L,\phi)$ denote the set of these elements in
$H^{2m}(K^m,K^m_{L})$, where $M$ is to vary over all thickenings
for which the $(m-\!1)$-st obstruction is zero. 
If there is an embedding of $K$ into
$B^4$, extending $\phi$, clearly there is a $4$-dimensional  
thickening $M$ which embeds into $S^4$, so 
$0\in o_m(K,L,\phi)$ for each $m$.

Consider the simplest relative case, when 
$(K,L)=(D^2\amalg\ldots\amalg D^2, S^1\amalg\ldots\amalg S^1)$.
By the result of Turaev \cite{Turaev}, the first non-trivial
obstruction coincides in this case with the first non-trivial
Milnor's $\bar\mu$-invariants of the link $\phi(L)$ in $S^3$.     
In this sense the obstructions $o_m(K)$ may be thought of
as an absolute analogue of $\bar\mu$-invariants. However,
since $2$-complexes in general have a more complicated topology,
$\{o_m(K)\}$ have a larger indeterminacy.

\section{Examples} \label{examples}

First we recall the construction of examples in \cite{FKT}.
Let $C$ denote the $2$-skeleton of the $6$-simplex with vertices 
$v_1,\ldots, v_7$, with one $2$-cell, with vertices $v_1 v_2 v_3$, 
removed. Take another copy $C'$ of $C$, with vertices $v'_1,\ldots,v'_7$, 
and denote by $\overline C$ the union of $C$ and $C'$, identified 
along their last vertices, $v_7=v'_7$. This $2$-complex is easily 
seen to admit an embedding into ${\mathbb{R}}^{4}$ (see \cite{Van Kampen}). 
Let $\gamma$ (resp. $\gamma'$) denote the loop $v_7 v_1 v_2 v_3 v_1 v_7$ 
(resp. $v_7 v'_1 v'_2 v'_3 v'_1 v_7$) in $\overline C$.

Denote by $F$ the free group on two generators, and
fix a positive integer $m$. Let $\alpha$ be an element in $F^m$, 
the $m$-th term of the lower central series of $F$. 
We identify $F$ with ${\pi}_1({\gamma}\vee{\gamma'})$, and we
associate to each word in $F$ its ``standard'' representative
loop in the wedge of two circles ${\gamma}\vee{\gamma'}$.
Finally, we construct the $2$-complex $K_{\alpha}$ by attaching
a $2$-cell to $\overline C$ along $\alpha$.

\begin{theorem}[\cite{FKT}] \label{FKTtheorem} \sl
Let $\alpha$ be a non-trivial element in $F^m$ for some $m\geq 2$. 
Then Van Kampen's obstruction $o(K_{\alpha})$ vanishes, but the $2$-complex 
$K_{\alpha}$ does not admit an embedding into ${\mathbb{R}}^4$.
\end{theorem}

We now present a computation of the obstructions $\{o_i(K_{\alpha})\}$.
The class $m$ of the commutator $\alpha$ is reflected in non-vanishing
of the obstruction $o_{m+1}(K_{\alpha})$.

\begin{theorem} \label{not-embeddable}  \sl
Let $\alpha$ be an element in $F^m$ for some $m\geq 2$,
and assume $\alpha\notin F^{m+1}$. Then $o_{m+1}(K_{\alpha})$ 
is defined and does not contain zero. 
In particular, $K_{\alpha}$ does not admit an embedding into ${\mathbb{R}}^4$.
\end{theorem}

\noindent
{\em Proof.}
First we construct a thickening $M$ of $K_{\alpha}$ with trivial 
intersection pairing and such that $o_i(K,M)=0$ for all $i\leq m$.
The complex $K_{\alpha}$ is obtained from $\overline{C}=C\vee C'$
by attaching a $2$-cell along the commutator $\alpha\in F^m$.
Van Kampen constructed in \cite{Van Kampen} an immersion of the
$2$-skeleton of the $6$-simplex with vertices $v_1,\ldots,v_7$ into
${\mathbb{R}}^4$ such that the $2$-cells with vertices $v_1 v_2 v_3$ and
$v_4 v_5 v_6$ intersect in one point, and all other simplices are 
disjoint and embedded. Consider the corresponding embedding of 
$\overline C$, and let $\overline M$ denote its thickening in
${\mathbb{R}}^4$. Clearly the intersection pairing on $\overline M$
is trivial, and all Massey products on $H^1(\partial \overline M)$
vanish. Recall that $C$ and $C'$ have the vertex $v_7$ in common.
Consider the handle decomposition of $\overline M$, given by
thickenings of simplices of $K_{\alpha}$ in ${\mathbb{R}}^4$.
The union of the handles in $\overline M$ corresponding to all simplices
in $\overline C$, containing $v_7$, is a $4$-ball $B$. The remaining 
$2$-handles are attached to $B$ along a link $\overline L$ in 
$S^3=\partial B$. Each attaching curve is isotopic, within the union of 
$0$- and $1$-handles, to the boundary curve of the corresponding
$2$-simplex. There is no $2$-cell attached to $v_1 v_2 v_3$,
however we introduce in $S^3$ a circle, isotopic to it.
Because of the choice of the embedding of $K$ into ${\mathbb{R}}^4$,
$\overline L$ is a slice link, and the curves isotopic to $v_1 v_2 v_3$
and $v_4 v_5 v_6$ (respectively $v'_1 v'_2 v'_3$ and $v'_4 v'_5 v'_6$)
have linking number one. The remaining $2$-cell of $K_{\alpha}$
is attached along the commutator of $v_1 v_2 v_3$ and $v'_1 v'_2 v'_3$.
Choosing appropriately the corresponding curve $l$ in $S^3$, we get
the link $L=\overline L \cup l$ such that all Milnor's $\bar\mu$-invariants
of $L$ of length less than $m+1$ vanish, and a $\bar\mu$-invariant
of length $m+1$ of the $3$-component link $(v_4 v_5 v_6, v'_4 v'_5 v'_6,
l)$ is non-trivial. It is a result of Turaev \cite{Turaev} that the first
non-vanishing $\bar\mu$-invariant is equal to the corresponding 
Massey product on $\partial M$, where $M=\overline M\cup_{l}2$-handle. 
(Note that the framings of the components
of $\overline L$ are zero, since the intersection pairing on $\overline M$
vanishes, and we choose the framing of $l$ also to be zero.)
This proves that $o_i(K,M)=0$ for all $i\leq m$, and 
$o_{m+1}(K_{\alpha},M)\neq 0$.

It remains to show that $o_{m+1}(K_{\alpha})$ does not contain zero.
Let $M$ be any thickening with trivial intersection pairing
and with $o_m(K,M)=0$. As above, the union of the handles in $M$ 
corresponding to all simplices, containing $v_7$, is a $4$-ball $B$.
By a theorem of Conway - Gordon \cite{CG}
and Sachs \cite{Sachs} any embedding of the complete graph on $6$ 
vertices in $S^3$ contains two disjoint linking cycles. Consider 
the complete graph on vertices $v_1,\ldots,v_6$ in $\overline C$.
According to the definition of $M$, the attaching curves of the
$2$-handles in $S^3=\partial B$ are isotopic to the attaching maps of 
simplices of $K_{\alpha}$. As above, we introduce in $S^3$ a curve 
isotopic to $v_1 v_2 v_3$.
Now we have in $S^3$ a perturbed version of the complete
graph on $6$ vertices. Since, according to definition \ref{thickening},
these perturbations take place in the union of $0$- and $1$-handles,
at least two of the curves must have a non-trivial linking number. 
Since the intersection pairing on $M$ vanishes, these two circles are 
the ones isotopic to $v_1 v_2 v_3$ and to $v_4 v_5 v_6$. Similarly we 
have in $S^3$ another, disjoint, copy of a perturbed graph on 
$v'_1,\ldots,v'_6$, and two linking circles isotopic to $v'_1 v'_2 v'_3$ 
and to  $v'_4 v'_5 v'_6$. Recall that there are no $2$-handles attached to
$v_1 v_2 v_3$ or $v'_1 v'_2 v'_3$, however there is a $2$-handle
whose attaching curve $l$ is a commutator of these circles.
It is easily seen that the link $(v_4 v_5 v_6, v'_4 v'_5 v'_6, l)$
has a non-trivial $\bar\mu$-invariant of length $m+1$. As above,
this is translated into non-vanishing of $o_{m+1}(K,M)$.
\qed

\begin{remark}
The idea of the proof of the fact that any embedding of the 
complete graph on $6$ vertices in $S^3$ contains two linking cycles 
(\cite{CG}, \cite{Sachs}) is conceptually similar to the proof of Van 
Kampen that the $2$-skeleton of the $6$-simplex does not embed into 
${\mathbb{R}}^4$ \cite{Van Kampen}. In both cases one shows that a certain
number is invariant mod $2$ for different maps - in one case,
the total linking number, in the other case, the total number of 
singular points of an immersion.
In this sense our proof of theorem \ref{not-embeddable} 
is similar to the proof of theorem \ref{FKTtheorem} in \cite{FKT}.
\end{remark}

\section{A note on configuration spaces} \label{config} 

In this section we consider an approach to the embedding problem,
suggested by obstruction theory and configuration spaces. We give a
reformulation of Van Kampen's obstruction in this context, which
suggests another approach to defining higher embedding obstructions.
Given a space $X$, $C^m(X)$ will denote its configuration space of $m$
points:

\[ C^m(X)=\{(x_1,\ldots,x_m)\in X^m\,|\,x_i\neq 
x_j\;\,\mathrm{if}\;\,i\neq j\}. \]

\noindent
In the simplicial category, for a complex $K$ we define

\[ C^m(K)=\{\sigma_1\times\ldots\times\sigma_m\subset
K^m\,|\,\mathrm{simplices}\;\,\sigma_i,\,\sigma_j \]
\[ \mathrm{have\;\,no\;\,vertices\;\,in\;\,common\;\,for}\;\,i\neq j\}. \]

\noindent
The configuration space of two points $C^2(X)$ is sometimes called deleted 
product and is also denoted by $X^{*}$.
The symmetric groups are denoted by $S_m$; $S_m$ acts freely
on $C^m(K)$, and on its $i$-skeleton $(C^m(K))^i$ for each $i$, 
by exchanging the coordinates.

A necessary condition for the existence of an embedding 
$K^n\hookrightarrow {{\mathbb{R}}}^{2n}$ is the existence, for each $m$,
of an $S_m$-equivariant map $C^m(K)\longrightarrow C^m({{\mathbb{R}}}^{2n})$.
We will now analyze the first embedding obstruction, corresponding to 
$m=2$, that is, the obstruction to existence of a 
${{\mathbb{Z}}}/2$-equivariant map 

\[ K\times K\setminus\Delta\longrightarrow 
{{\mathbb{R}}}^{2n}\times{{\mathbb{R}}}^{2n}\setminus\Delta\simeq S^{2n-1}.\]

\noindent
The ${\mathbb{Z}}/2$-equivariant homotopy equivalence above is given by
the projection of ${\mathbb{R}}^{2n}\times {\mathbb{R}}^{2n}\setminus\Delta$
onto the unit sphere in the antidiagonal $\{(x,-x)\}\subset 
{\mathbb{R}}^{2n}\times {\mathbb{R}}^{2n}$.
The diagonal $\Delta$ in $K\times K$ is the ``simplicial'' diagonal,
as defined in section \ref{VKobstruction}, while $\Delta\subset
{\mathbb{R}}^{2n}\times {\mathbb{R}}^{2n}$ is the usual set-theoretic
diagonal. Recall that the spaces above are denoted in short by
$K^*$ and $({\mathbb{R}}^{2n})^*$ respectively.  

\begin{theorem} \label{reformulation1}  \sl
The obstruction to existence of a ${\mathbb{Z}}_2$-equivariant map
$K^{*}\longrightarrow ({\mathbb{R}}^{2n})^{*}$ lies in 
$H^{2n}_{{{\mathbb{Z}}}/2}(K^*;{{\mathbb{Z}}})$
and coincides with Van Kampen's obstruction $o(K)$. 
\end{theorem}

\noindent
{\em Proof.} 
Since $K^{*}$ is a $(2n)$-dimensional CW-complex,
the only non-trivial obstruction group in this setting is
$H^{2n}_{{\mathbb{Z}}/2}(K^{*};$ $\pi_{2n-1}(S^{2n-1}))\cong
H^{2n}_{{{\mathbb{Z}}}/2}(K^*;{{\mathbb{Z}}})$.

Let $f: K\longrightarrow {\mathbb{R}}^{2n}$ be any immersion.
Since $f(\sigma)$ and $f(\nu)$ are disjoint for any $n$-simplex
$\sigma$ and any $(n-1)$-simplex $\nu$, $f\times f$ restricted to the
$(2n-1)$-skeleton of $K^*$   
is a ${\mathbb{Z}}/2$-equivariant embedding into 
$({\mathbb{R}}^{2n})^*$.
Let $\sigma$, $\tau$ be two $n$-dimensional simplices of $K$ and consider
$\sigma\times\tau$ as an oriented generator of $(2n)$-dimensional
cellular chains on $K^{*}$. The obstruction cochain $c_f$ assigns to
$\sigma\times\tau$ the element 

\[ c_f(\sigma\times\tau)=
[(f\times f)(\partial(\sigma\times\tau))]\in \pi_{2n-1}(S^{2n-1}). \]

\noindent 
The map $f\times f$ sends $\sigma\times\tau$ into 
${\mathbb{R}}^{2n}\times{\mathbb{R}}^{2n}$, and one has

\[ o_f(\sigma\times\tau)=f(\sigma)\cdot f(\tau)=
(f\times f)(\sigma\times\tau)\cap\Delta_{{\mathbb{R}}^{2n}}=(f\times
f)(\partial(\sigma\times\tau))=c_f(\sigma\times\tau). \]

\noindent
This shows that the homotopy-theoretic obstruction coincides
with Van Kampen's obstruction even on the cochain level,
when the map of $(2n-1)$-skeleton of $K^*$ corresponds
to the chosen immersion $f$.
This completes the proof, since the cohomology class $[c_f]$ is 
independent of the 
choice of a map of $(2n-1)$-skeleton of $K^*$, being the first non-trivial
obstruction.
\qed

\begin{remark}
This result is implicitely contained in \cite{Haefliger}, \cite{Shapiro}, \cite{Wu}. 
It is interesting to note that by theorems \ref{iff} and
\ref{reformulation1}, the existence of a ${\mathbb{Z}}/2$ - equivariant map
$K^{*}\longrightarrow ({\mathbb{R}}^{2n})^{*}$ is equivalent
to existence of an embedding $K^{n}\hookrightarrow {{\mathbb{R}}}^{2n}$
for $n\neq 2$.
\end{remark}

We conclude by suggesting the following approach to defining higher
embedding obstructions, which will be pursued in a separate paper. 
Suppose $o(K)$ vanishes, so there exists
a ${\mathbb{Z}}/2$-equivariant map $K^*\longrightarrow ({\mathbb{R}}^4)^*$.
One may consider the obstructions to existence of an equivariant,
with respect to the free action of the symmetric group $S_m$,
map $C^m(K)\longrightarrow C^m({\mathbb{R}}^4)$, $m=3,4,\ldots$. 
Fadell and Neuwirth \cite{FN} have determined homotopy types 
of the symmetric products of Euclidean spaces, thus (rationally) 
explicitely giving the coefficients of obstruction groups. Note that
the examples in \cite{SSS} (similar to those constructed in \cite{FKT}) 
show that the entire sequence of such obstructions, arising from 
configuration spaces, is incomplete. 

{\em Acknowledgements.} I would like to thank Michael Freedman, Richard
Stong and Peter Teichner for many discussions.

This work has been done during my stays at the University of
California - San Diego, Michigan State University and the Max-Planck-Institut
f\"{u}r Mathematik, and I would like to thank them for their hospitality
and support.

\end{document}